\documentclass{amsart}
\usepackage{amsmath,amssymb,amscd}
\usepackage{color}
\usepackage[varg]{txfonts}

\newtheorem{theorem}{Theorem}[section]

\newtheorem{corollary}[theorem]{Corollary}

\newtheorem{fact}[theorem]{Fact}
\newcommand{\COMP}{\raisebox{0.3ex}{\hspace{0.3ex}%
 {$\scriptstyle{\circ}$}\hspace{0.5ex}}}   
\newcommand{\FRAC}[2]{\displaystyle\frac{#1}{#2}}
\newcommand{\OL}[1]{\overline{#1}}
\renewcommand{\phi}{\varphi}

\newcommand{\trdeg}{\operatorname{trdeg\,}}
\newcommand{\vv}{\boldsymbol{v}}

\newcommand{\vvv}{\hat{\boldsymbol{v}}}
\title[Diophantine inequality for local domains]
{Diophantine inequality for equicharacteristic\\
excellent Henselian local domains}
\author{Hirotada ITO, \ \ Shuzo IZUMI}
\date{\today}
\begin{document}
\maketitle
\begin{center}
\begin{minipage}{65ex}
\noindent
\textit{A{\scriptsize BSTRACT}.} 
G. Rond has proved a Diophantine type inequality for the field of 
quotients of the convergent 
or formal power series ring in multivariables. 
We generalize his theorem to the field of the quotients of an excellent 
Henselian local domain in equicharacteristic \vspace{2ex}case. 

\noindent
\textit{R{\scriptsize\'{E}SUM\'{E}}.} 
G. Rond a d\'{e}montr\'{e} une in\'{e}galit\'{e} de type Diophantien pour le 
corps des 
quotients de s\'{e}ries convergentes (ou formelles) \`{a} plusieurs 
variables. 
On fait ici une g\'{e}n\'{e}ralisation de son th\'{e}or\`{e}me au corps des 
quotients d'un anneau local int\'{e}gral henselien excellent 
dans le cas \'{e}qui-caract\'{e}ristique.
\vspace{.3cm} 
\end{minipage}

 \vspace{1ex}13F40, 
\vspace{1ex}13H99, \vspace{1ex}11D75, \vspace{1ex}11J25

K{\scriptsize EY WORDS}: 
Diophantine inequality, m-valuation,\\
associated valuations, linear Artin approximation property. 
\end{center}
\section{Introduction}
An important topic of Diophantine approximation is the 
problem of approximation of a real algebraic number by rational ones. 
The crucial result is Roth's theorem:

\textit{
If $z\in\mathbb{R}\setminus\mathbb{Q}$ is an algebraic number, 
}
$$
\forall\epsilon>0\ 
\exists c(z,\ \epsilon)>0\ \forall x\in\mathbb{Z}\ \forall y\in\mathbb{Z^*}:\ 
\bigl|z-\frac{x}{y}\bigr|>c(z,\ \epsilon)|y|^{-2-\epsilon}.
$$
There are quite similar results for the Laurent series field 
in a single variable (cf. \cite{lasjaunias}). 
It is also known that there are deep analogous results 
on the global function fields on certain special varieties 
in connection with Nevanlinna's theory (cf. \cite{ru}). 

Rond \cite{rond2} obtained a Diophantine inequality for the field 
of quotients of the convergent 
or formal power series ring in multivariables in connection with the 
\textit{linear Artin approximation property} 
(Spivakovsky, cf. \cite{rond1}). 
He used the \textit{product inequality} \cite{izumi1} 
for the order function $\nu$ on an analytic integral domain. 

In this paper we assert that Diophantine inequality holds for the field of 
quotients of an equicharacteristic excellent Henselian local domain. 
For the proof, we need \textit{Rees's inequality \cite{rees4} 
for $\mathfrak{m}$-valuations} on complete local rings, 
a variant of the product inequality. 
To be precise, we use its generalization 
to analytically irreducible 
excellent domains by H\"{u}bl-Swanson \cite{hs}. 

An inequality on the order function was once used 
for zero-estimate of elements transcendental over the polynomial ring 
generated by a parameter system in a local ring in \cite{izumi2}. 
This time we are concerned with elements algebraic over a local ring. 

Let us give a precise description of our theorem. 
Let $K$ be a (commutative) field. We call a mapping 
$\vv :K\longrightarrow \OL{\mathbb{R}}\ 
(\OL{\mathbb{R}}:=\mathbb{R}\cup\{\infty\})$ 
a \textit{valuation} when it satisfies the following:
\begin{enumerate}
\item
$\vv (xy)=\vv (x)+\vv (y)$,
\item
$\vv (x+y)\ge\min\{\vv (x),\ \vv (y)\}$,
\item
$\vv(0)=\infty$.
\end{enumerate}
We can define the absolute value 
$|\ |_{\vv}:K\longrightarrow \mathbb{R}$ by 
$|w|_{\vv}:=\exp(-{\vv }(w))$ $(|0|_{\vv}:=0)$. 
Then $K$ is a metric space defined by the absolute value of 
the difference. This defines a topology compatible with the field operations.
We endow the discrete topology on $\mathbb{R}$. 
Then $\vv:K^*\longrightarrow\mathbb{R}$ $(K^*:=K\setminus 0)$ is continuous. 
We put:
\\ \hspace*{3ex}
$V_{\vv}:=\{z\in K:\vv(z)\ge 0\}$: the \textit{valuation ring} 
of $\vv $,
\\ \hspace*{3ex}
$\mathfrak{m}_{\vv}:=\{z\in K:\vv(z)>0\}$: the \textit{maximal 
ideal} of $V_{\vv}$,
\\ \hspace*{3ex}
$k_{\vv}:=V_{\vv}/\mathfrak{m}_{\vv}$: 
the \textit{residue field} of $\vv$,
\\ \hspace*{3ex}
$\hat{K}$: the completion of $K$ with respect to $\vv$, 
which has a natural structure of a field,
\\ \hspace*{3ex}
$\vvv$: the continuous extension of $\vv$ to $\hat{K}$, which is 
a valuation on $\hat{K}$, 
\\ \hspace*{3ex}
$\hat{V}_{\vv}$: the valuation ring of the extension $\vvv$. 

A valuation is called \textit{discrete valuation} if the value group 
$\vv(K^*)$ is isomorphic to $\mathbb{Z}$ as an ordered group. 
In this case 
the valuation ring $V_{\vv}$ is a discrete valuation ring (DVR) and 
$\vv $ coincides with the $\mathfrak{m}_{\vv}$-adic order on $V_{\vv}$ 
and we have $K=Q(V_{\vv})$. 
The completion $\hat{B}$ of some subset $B\subset K$ can be identified with 
its closure in $\hat{K}$. 
Our main result is the following. 
\\

\textit{
Let $(A,\ \mathfrak{m})$ be an 
equicharacteristic analytically irreducible excellent Henselian local 
domain and $\vv$ an $\mathfrak{m}$-valuation (defined in the next section) 
on the field $K:=Q(A)$ of quotients of $A$. 
If $z\in\hat{K}\setminus K$ is algebraic over $K$, 
then we have the following:} 
$$
\exists a > 0\ \exists c> 0\ \forall x\in A\ \forall y\in A^*:\ 
\bigl|z-\frac{x}{y}\bigr|_{\vvv}>c\,|y|_{\vv }^a
.$$

Note that $\hat{K}$ is not generally the field quotients of $\hat{A}$ 
(cf. \cite{rond1}, {\bf 2.4}). 
The essential point of the proof is reducing inequality on the 
valuation to inequality on the maximal-ideal-adic order in the same way 
as in \cite{rond2} (see \S 4, (v)). 
In our general case, 
we need Rees's valuation theorem \cite{rees2}, 
\cite{rees3} to connect valuations to the order. 

In the contrary to the case of 
algebraic numbers, the exponent on the right of this inequality 
is not uniformly bounded. Rond (\cite{rond1}, {\bf 2.4}) has 
given a sequence of elements $z_i\in\hat{K}$ 
of degree 2 over $K$ with unbounded exponents. 
\section{$\mathfrak{m}$-valuations on local domains}\label{mvaluation}
Let $(A,\ \mathfrak{m})$ be a local domain 
whose field of quotients $Q(A)$ is $K$. 
Let $k:=A/\mathfrak{m}$ denote the residue field. 
A valuation $\vv $ on $K$ is called an $\mathfrak{m}$-\textit{valuation}, 
if it satisfies the following:

(a) $x\in A\Longrightarrow\vv (x)\ge 0$,

(b) $x\in\mathfrak{m}\Longrightarrow\vv (x)>0$,

(c) $\trdeg_k k_{\vv}=\dim A-1$,

(d) The value group $\vv(K^*)$ is isomorphic to $\mathbb{Z}$ 
(as an ordered group). 
\vspace{1ex}

Let us recall key facts on valuations which is used in the proof. 
The first one is Rees's 
\textit{strong valuation theorem} \cite{rees3}. 
We state only the special case which we need later. 
We define the $\mathfrak{m}$\textit{-adic order} 
$\nu_{\mathfrak{m}} :A\longrightarrow\OL{\mathbb{R}}$ on $A$ 
by $\nu_{\mathfrak{m}} (f):=\max\{p:f\in\mathfrak{m}^p\}$. 
This is not necessarily a valuation. It satisfies formulae 
\\ \hspace*{2ex}
$(1')$\  $\nu_{\mathfrak{m}} (fg)\ge\nu_{\mathfrak{m}} (f)
+\nu_{\mathfrak{m}} (g)$,
\\ \hspace*{2ex}
$(2)$\  $\nu_{\mathfrak{m}} (f+g)\ge
\min\{\nu_{\mathfrak{m}} (f),\ \nu_{\mathfrak{m}} (g)\}$,
\\ \hspace*{2ex}
$(3')$\ $\nu_{\mathfrak{m}}(0)=\infty$, $\nu_{\mathfrak{m}}(1)=0$. 

Let us stabilize $\nu_{\mathfrak{m}}$ by Samuel's idea: 
$\OL{\nu}_{\mathfrak{m}}(f):=\lim_{k\to\infty}\nu_{\mathfrak{m}}(f^k)/k.$ 
This limit always exists and 
satisfies formulae $(1')$, $(2)$, $(3')$ and the homogeneity formula
\\ \hspace*{2ex}
(4)\ $\OL{\nu}_{\mathfrak{m}}(f^n)=n\,\OL{\nu}_{\mathfrak{m}}(f)
\quad (n\in\mathbb{N})$\\
also (see \cite{rees1}). The following is Rees's 
\textit{Strong valuation theorem}.  
\begin{fact}[\cite{rees2}, \cite{rees3}]\label{valuation}
Let $(A,\ \mathfrak{m})$ be a Noetherian local ring 
whose $\mathfrak{m}$-adic completion is reduced 
(has no non-zero nilpotent element). 
Then there exist a non-negative number $C$ 
and a set of valuations $\vv_1,\dots,\vv_p$ on $K$ with the value group 
$\mathbb{Z}$ such that 
$$
\forall x\in A:\ 
\nu_{\mathfrak{m}}(x) 
\le\OL{\nu}_{\mathfrak{m}}(x)
\le\nu_{\mathfrak{m}}(x)+C,
$$$$
\forall x\in A:
\OL{\nu}_{\mathfrak{m}}(x)
=\min\{r_1\vv_1(x),\dots,r_p\vv_p(x)\}\ \ 
(r_i:=1/\min\{\vv_i(y):y\in\mathfrak{m}\}).
$$
The set $\{\vv_1(x),\dots,\vv_p(x)\}$ is unique, 
if it is taken irredundant.
\end{fact}
We call the irredundant valuations $\vv_1,\dots,\vv_p$ the 
\textit{valuations associated with} $\mathfrak{m}$. 
We call a local ring \textit{analytically irreducible} when 
its $\mathfrak{m}$-adic completion is an integral domain. 
Rees proves the following: 

\begin{fact}[\cite{rees1}, {\bf 5.9}]\label{equicharacteristic}
Let $(A,\ \mathfrak{m})$ be an equicharacteristic 
analytically irreducible local domain. 
Then the valuations associated with $\mathfrak{m}$ are all 
$\mathfrak{m}$-valuations. 
\end{fact}
In the proof of the regular analytic case, Rond \cite{rond2} use the product 
inequality \cite{izumi1} for analytic domain. 
Rees generalises this inequality and, in the complete domain case, 
gives a valuation theoretic form \cite{rees4}, 
({\bf E}). H\"{u}bl and Swanson generalise the latter to excellent domains 
as follows: 
\begin{fact}[\cite{hs}, {\bf 1.3}]\label{rees}
Let $(A,\ \mathfrak{m})$ be an 
analytically irreducible excellent local domain. 
Then for any pair of $\mathfrak{m}$-valuations $\vv$ and $\vv'$, 
we have the following.
$$
\exists d > 0\ \forall x\in A:\ 
\vv(x)\le d\,\vv'(x).
$$
The constant $d$ can be chosen independent of $\boldsymbol{v}'$.
\end{fact}
Combining these facts we see the following: 
\begin{fact}\label{comparison}
Let $(A,\ \mathfrak{m})$ be an equicharacteristic 
analytically irreducible excellent local domain and let $\vv$ 
be an $\mathfrak{m}$-valuation on $A$. Then we have:
$$
\exists C>0\ \exists s>0\ \exists t>0\ \forall x\in A:\ 
s\,\vv(x)
\le\OL{\nu}_{\mathfrak{m}}(x)
\le{\nu}_{\mathfrak{m}}(x)+C
\le t\,\vv(x)+C
.$$
\end{fact}

\section{Main theorem}
With the notation in Introduction our main theorem is the following. 
\begin{theorem}\label{thm}
Let $(A,\ \mathfrak{m})$ be an equicharacteristic analytically 
irreducible excellent Hen\-selian local domain 
and let $K:=Q(A)$ denote its field of quotients and let 
$\vv :K\longrightarrow\OL{\mathbb{R}}$ 
be an $\mathfrak{m}$-valuation. 
If $z\in \hat{K}\setminus K$ is algebraic over $K$, then we have
$$ 
\exists a> 0\ \exists c> 0\ \forall x\in A\ \forall y\in A^*:\ 
\bigl|z-\frac{x}{y}\bigr|_{\vvv }>c|y|_{\vv }^a
.$$
\end{theorem}
Just in the same way as Rond \cite{rond2}, {\bf 3.1} 
(see \cite{rond1}, {\bf 2.1} also), 
our Theorem \ref{thm} implies the following. 
\begin{corollary}\label{cor}
Let $(A,\ \mathfrak{m})$ be an equicharacteristic analytically irreducible 
excellent Hen\-selian domain and let 
$P(X,\ Y)\in A[X,\ Y]$ be a homogeneous polynomial. 
Then the Artin function of $P(X,\ Y)$ is majorised by an affine 
function, i.e. 
$$
\exists\alpha\ \exists\beta\ \forall x\in A\ \forall y\in A:\ 
\nu_{\mathfrak{m}}(P(x,\ y))\ge\alpha i + \beta
$$$$
\Longrightarrow
\exists\OL{x}\in A\ \exists\OL{y} \in A:\ 
\nu_{\mathfrak{m}}(\OL{x}-x)\ge i,\ \nu_{\mathfrak{m}}(\OL{y}-y)\ge i,\ 
P(\OL{x},\ \OL{y})=0
.$$
\end{corollary}
This corollary reminds us of the theorem that an excellent Henselian local 
ring has the strong Artin approximation property (cf. \cite{popescu}). 
The case $P(X,\ Y)=XY$ is nothing but the product inequality. 
\section{Proof of Theorem}
\textbf{(i) Reduction to normal case.} 

We may assume that $\vv(K^*)=\vv(\hat{K}^*)=\mathbb{Z}$. 
This results in a change of the exponent $a$. 
Let $\tilde{A}$ denote the normalization (the integral closure 
of $A$ in $K$) of $A$. 
Since $A$ is a Henselian integral domain, $\tilde{A}$ is a local 
ring by \cite{nagata}, {\bf 43.11}, {\bf 43.20}. 
Since $A$ is excellent, 
\begin{equation}
A\text{ is a G-ring and a Nagata ($=$ pseudo-geometric) ring}
\end{equation} 
by \cite{matsumura}, {\bf 33.H}. 
Then $\tilde{A}$ is a finite $A$-module. Hence 
$\dim A=\dim\tilde{A}$ 
by a theorem of Cohen-Seidenberg (cf. \cite{nagata}, {\bf 10.10}) 
and $r\tilde{A}\subset A$ 
for some $r\in A^*$ (existence of a universal denominator). 
Then a Diophantine inequality for $\tilde{A}$ implies one for $A$ 
with the same exponent $a$. 
Finiteness also implies that $\tilde{A}$ is excellent 
and Henselian by \cite{nagata}, {\bf 43.16}. 

Let $\tilde{\mathfrak{m}}$ denote the maximal ideal of $\tilde{A}$. 
We claim that $\vv$ is an $\tilde{m}$-valuation. 
If $x\in\tilde{A}$, 
$$
\exists p\in\mathbb{N},\ \exists b_0,\dots,b_{p-1}\in A:\ 
x^p=b_o+b_1x+\dots+b_{p-1}x^{p-1}
.$$
Then we have 
$p\vv(x)\ge\min\{ i\vv(x):0\le i\le p-1 \}.$ 
This proves $\tilde{A}\subset V_{\vv}$ and condition (a) for 
$(\tilde{A},\ \tilde{\mathfrak{m}})$. Let us put 
$\OL{\mathfrak{m}}:=\{x\in\tilde{A}:\vv(x)>0\}.$ 
Then $\OL{\mathfrak{m}}$ is a prime ideal of $\tilde{A}$ and 
$
\OL{\mathfrak{m}}\cap A
=\mathfrak{m}.
$ 
This implies that 
$\OL{\mathfrak{m}}=\tilde{\mathfrak{m}}$ by \cite{AC}, Chapt.\/5, 2.1, 
Prop.~1 and (b) holds. 
Since $\tilde{A}$ is a finite $A$-module, 
$\tilde{k}=\tilde{A}/\tilde{\mathfrak{m}}$ is a finite
$k$-module $(k:=A/\mathfrak{m})$ \textit{i.e.} $\tilde{k}$ is algebraic over $k$. 
This proves (c). The condition (d) is obvious. 
We have proved the claim and we may assume that 
\begin{equation}
A\text{ is an equicharacteristic, excellent, Henselian and 
normal local domain.}
\end{equation}

\textbf{(ii) Reduction of the minimal equation.}

Let
$$
\phi(Z):=a_0+a_1Z+\dots+a_dZ^d\  (a_d\neq 0,\ d\ge 2)
$$
be a minimal equation for $z$ over $A$, that is, $\phi$ 
is a polynomial of the minimal degree in $A[Z]$ with $\phi(z)=0$. 
Now take $u\in A^*$ and put 
$$
\phi_u(Z):=u^da_d^{d-1}\phi({Z}/{ua_d}).
$$ 
Then we have 
$$
\phi_u(Z)=a_0u^da_d^{d-1}+a_1u^{d-1}a_d^{d-2}Z+\dots+Z^d\in A[Z]
$$
and $w'\in \hat{K}$ is a root of $\phi_u(Z)$ if and only if 
$w:={w'}/{ua_d}$ is 
a root of $\phi(Z)$. 
If 
$$
\exists a\ge 0\ \exists c> 0\ \forall x\in A\ \forall y\in A:\ 
\bigl|z'-\FRAC{x}{y}\bigr|_{\vvv}>c|y|_{\vv}^a
$$ 
holds for $z':=ua_dz$, we have
$$
\exists a\ge 0\ \exists c> 0\ \forall x\in A\ \forall y\in A:\ 
\bigl|z-\FRAC{x}{y}\bigr|_{\vvv}>\FRAC{c}{|ua_d|}|y|_{\vv}^a
.$$ 
The polynomial $\phi_u(Z)\in A[Z]$ is minimal for $z'$. 
Thus, choosing $u$, we may assume that $z\in \hat{V}_{\vv}$ and
\begin{equation}
\phi(Z):=a_0+a_1Z+\dots+a_{d-1}Z^{d-1}+Z^d\  
(d\ge 2,\ a_i\in\mathfrak{m}^{d-i})\label{integral}
\end{equation}
from the first. 

\textbf{(iii) Order function on $A[z]$ 
(the ring generated by $z$ over $A$).}

Let us consider the residue ring $B:=A[Z]/{\phi(Z)A[Z]}$. 
There is an isomorphism $\iota:B\longrightarrow A[z]$. 
The ring $B$ is a finite $A$-module with basis $1,\ z,\ z^2,\dots,\ z^{d-1}$. 
Since $A[z]$ is a subring of the field $\hat{K}$, $B$ is an integral 
domain. Thus we have the following:
\begin{equation}
A[z]\cong B:=A[Z]/\phi(Z) A[Z]\text{ is an integral extension of $A$.} 
\label{integralext}
\end{equation}

Since $A$ is Henselian, $B$ is a local ring by (\ref{integralext}) 
and by \cite{nagata}, {\bf 43.12}. As a consequence of 
(ii), $z^d\in\mathfrak{m}A[z]$. 
Hence the maximal ideal of $B$ is $\mathfrak{n}:=\mathfrak{m}B+ZB$ 
and its residue ring is the same as that of $A$: 
$k=A/\mathfrak{m}=B/\mathfrak{n}$. 
Let us define $\mu:A[Z]\longrightarrow \OL{\mathbb{R}}$ by
$$
\mu(\sum_{i=0}^{e}b_iZ^i)
:=\underset{i}{\min}\{\nu_{\mathfrak{m}}(b_i)+i\}\ 
(b_i\in A)
$$
and $\nu_{\mathfrak{n}}:B\longrightarrow \OL{\mathbb{R}}$ 
as the $\mathfrak{n}$-adic order. The function $\mu$ is nothing but the 
restriction of the standard order on the formal power series ring $A[[Z]]$. 
We claim that $\nu_{\mathfrak{n}}(x)$ coincides with the 
$\mu$-order of the unique representative of $x$ 
in $A[Z]$ of degree less than $d$ \textit{i.e.} 
$$
\mu\bigl(\sum_{i=0}^{d-1} b_iZ^i\bigr)
=
\nu_{\mathfrak{n}}\bigl(\sum_{i=0}^{d-1} b_iZ^i\hspace{-.9ex}
\mod\phi(Z) A[Z]\bigr)
.$$
We have only to show that inequality 
$$
{\mu}\bigl(\sum_{i=0}^{d-1}b_iZ^i\bigr)
<{\mu}\bigl(\sum_{i=0}^{d-1}b_iZ^i
+\sum_0^da_iZ^i\sum_{j=0}^e c_jZ^j\bigr)
$$
leads us to a contradiction. Let us develop the product 
$\sum_{i=0}^da_iZ^i\sum_{i=0}^e c_jZ^j$ and reduce its degree in 
$Z$ by repeated substitutions $Z^d=-\sum_{i=0}^{d-1}a_iZ^i$, 
beginning from the highest degree term. 
By the assumption $a_i\in\mathfrak{m}_A^{d-i}$, 
the substitutions do not lower the $\mu$-order and we 
reach the left side. This contradicts the inequality we assumed.  

The function ${\nu}_{\mathfrak{n}}$ induces 
$\nu:={\nu}_{\mathfrak{n}}\COMP\iota^{-1} :
A[z]\longrightarrow \OL{\mathbb{R}}.$ 
Of course $\nu$ inherits the \textit{non-cancellation property} from 
${\nu}_{\mathfrak{n}}$: 
$$
\nu\bigl(\sum_{i=0}^{d-1} b_iz^i\bigr)
=
\min\{\nu(b_i)+i:0\le i\le d-1\}
=
\min\{\nu_{\mathfrak{m}}(b_i)+i:0\le i\le d-1\}
.\vspace{1ex}$$
In other words, there occurs no cancellation among summands 
of degree less than $d$. 

\textbf{(iv) $A[z]$ is analytically irreducible.}

Since $A$ is a normal G-ring, it is analytically normal, 
\textit{i.e.} the completion $\check{A}$ 
with respect to the $\mathfrak{m}$-adic topology is normal 
by \cite{matsumura}, {\bf 33.I}. Hence, by {\bf 2.4}, we have 
\begin{equation}
\check{A}=\hat{A}\text{ is normal} 
\label{analyticallynormal}.
\end{equation}
Let $\mathfrak{m}'$ denotes the maximal ideal of $A[z]$.
Taking the last equality of (iii) into account, 
we see that the $\mathfrak{m}'$-adic completion of $A[z]$ 
is isomorphic to $\hat{A}[z] \subset \hat{K}$. 
Hence $A[z]$ is an analytically irreducible domain. 
(This can be also deduced from \cite{nagata}, {\bf 44.1} and 
(1), (\ref{integralext}), (\ref{analyticallynormal}).) 
Now we can apply \ref{comparison} to $A[z]$. 
 
\textbf{(v) Diophantine inequality.}

We claim that the restriction $\vvv|_{Q(A[z])}$ is an 
$\mathfrak{m}'$-valuation. 
By the reduction (ii) 
we see that $z\in \hat{V}_{\vv}$ using the argument 
in (i) and (a) follows. 
Since $\mathfrak{m}'$ is generated by $\mathfrak{m}$ and $z$, 
the conditions (b) is satisfied. 
Take any element $x\in\hat{V}_{\vv}\cap K[z]$. 
There exists a nontrivial polynomial relation 
$$
c_0+c_1x+\dots+c_{p-1}x^{p-1}+c_px^p=0\quad(c_i\in A).
$$ 
If $t$ is a generator of $\mathfrak{m}_{\vv}\subset V_{\vv}$, 
we have an expression 
$$
c_i=c'_it^{\alpha_i}\quad 
(\alpha_i\in\{0,1,2,\dots\},\ 
c'_i\in V_{\vv}\setminus\mathfrak{m}_{\vv})
.$$ 
We may assume that some $\alpha_i$ is zero. 
Then the equation implies that 
$x\hspace{-.9ex}\mod\!\hat{\mathfrak{m}}_{\vv}\cap K[z]$ is algebraic 
over $k_{\vv}=V_{\vv}/\mathfrak{m}_{\vv}$. Therefore 
$(\hat{V}_{\vv}\cap K[z])/(\hat{\mathfrak{m}}_{\vv}\cap K[z])$ 
is algebraic over $k_{\vv}$ and we have 
$$
\trdeg_{A[z]/\mathfrak{m}'}
(\hat{V}_{\vv}\cap K[z])/(\hat{\mathfrak{m}}_{\vv}\cap K[z])
=\trdeg_k
(\hat{V}_{\vv}\cap K[z])/(\hat{\mathfrak{m}}_{\vv}\cap K[z])
$$
$$
=
\trdeg_k k_{\vv}=\dim A-1=\dim A[z]-1
.$$ 
Here the third equality follows from (c) for $\vv$. 
This equality implies (c) for $\vvv|_{K[z]}$. 
The condition (d) is obvious. Thus we have proved the claim. 

If ${\vvv}\bigl(z-\FRAC{x}{y}\bigr)\le \vvv(z)$, we have 
$\bigl|z-\FRAC{x}{y}\bigr|_{\vvv}\ge \exp(-\vvv(z))$
at once. 
Hence we may assume that ${\vvv}(x-yz)-\vv(y)>\vvv(z).$ 
If $\vv(x)\neq\vvv(yz)$, we have a contradiction: 
$$
\vvv(z)<\vvv(x-yz)-\vv(y)\le\vvv(yz)-\vv(y)=\vvv(z).
$$
Consequently we have only to consider the case
$$
\vv(x)=\vvv(yz)=\vv(y)+\vvv(z)
.$$
Since $A[z]$ is analytically irreducible, 
applying the inequality \ref{comparison} 
and the equality at the last part of (iii), we have
$$
s\,\vvv(x- yz)
\le
\nu_{\mathfrak{m}'}(x-yz)+C
\le
\nu_{\mathfrak{m}'}(x)+C
\le
t\,\vv(x)+C.
$$
It follows that 
$$
s\,\vvv\bigl(z-\FRAC{x}{y}\bigr)
\le
t\,\vv(x)-s\,\vv(y)+C
=
(t-s)\,\vv(y)+t\,\vvv(z)+C
.$$
This implies the inequality of our theorem. 

If $a=0$, $z$ is isolated from $K$ and cannot be in $\hat{K}$. 
Hence we see that $a>0$. 

\noindent
\begin{minipage}[t]{15.7em}
Interdisciplinary Graduate School\\
of Science and Engineering,\\
Kinki University\\
Higashi-Osaka 577-8502, Japan
\end{minipage}
\hspace{5em}
\begin{minipage}[t]{14em}
Department of Mathematics\\
Kinki University\\
Higashi-Osaka 577-8502, Japan
\end{minipage}
\end{document}